\documentclass[12pt,twoside]{article} 
\usepackage{amssymb,mathrsfs}
\usepackage{graphicx}
\usepackage[small,sc]{caption} 

\newtheorem{theorem}{Theorem}[section]
\newtheorem{lemma}[theorem]{Lemma}
\newtheorem{proposition}[theorem]{Proposition}
\newtheorem{corollary}[theorem]{Corollary}
\newtheorem{definition}[theorem]{Definition}

\newtheorem{rmrk}[theorem]{Remark}

\makeatletter
\@addtoreset{equation}{section}
\makeatother

\newcommand{\fig}[3] {
\medskip\smallskip
\begin{figure}[htb]
  \centering
  \includegraphics[width=#2]{#1.pdf}
  \begin{minipage}[t]{0.80\linewidth} 
    \caption{#3}
    \protect\label{#1}
  \end{minipage}
\end{figure}
\medskip
}

\newenvironment{remark}
{\begin{rmrk} \em}
{\end{rmrk}}

\newcommand{\fn} {function}
\newcommand{\me} {measure}
\newcommand{\tr} {trajector}
\newcommand{\erg} {ergodic}
\newcommand{\sy} {system}

\newcommand{\pr} {probability}
\newcommand{\ra} {random}
\newcommand{\dsy} {dynamical system}
\renewcommand{\o} {orbit}

\newcommand{\R} {\mathbb{R}}

\newcommand{\Z} {\mathbb{Z}}
\newcommand{\N} {\mathbb{N}}
\newcommand{\qed} {\hfill {\small Q.E.D.} \par\medskip}
\newcommand{\skippar} {\par\medskip}

\newcommand{\proof} {\noindent \textsc{Proof.} }
\newcommand{\proofof}[1] {\noindent \textsc{Proof of {#1}.} }
\newcommand{\article}[3] {\textsc{{#1}}, {\itshape {#2}}, {{#3}}.}
\newcommand{\book}[3] {\textsc{{#1}}, {\itshape {#2}}, {{#3}}.}
\newcommand{\vol} {\textbf}
\newcommand{\eps} {\varepsilon}
\newcommand{\rset}[2] {\left\{ #1 \: \left| \: #2 \right. \! \right\} }

\renewcommand{\iff} {if and only if\ }

\newcommand{\rw} {random walk}
\newcommand{\en} {environment}

\newcommand{\om} {\omega}
\newcommand{\ps} {\mathcal{M}}
\newcommand{\ma} {T}
\newcommand{\im} {\mu_*}
\newcommand{\ime} {\Pi_*}
\newcommand{\sime} {\Omega_*}
\newcommand{\ima} {\mathbb{P}_*}

\newcommand{\ld} {\mathcal{D}}
\newcommand{\gdd} {g.$d$-d.}
\newcommand{\Omd} {\Omega_\mathrm{d}}
\newcommand{\Omt} {\Omega_\mathrm{t}}

\begin{document}

\title{\textbf{Random walks in random environments without 
ellipticity}}

\author{\textsc{Marco Lenci}
\thanks{
Dipartimento di Matematica, 
Universit\`a di Bologna, 
Piazza di Porta S.\ Donato 5, 
40126 Bologna, Italy.
E-mail: \texttt{marco.lenci@unibo.it} } 
\thanks{INFN, Sezione di Bologna,
Via Irnerio 46,
40126 Bologna, Italy.}
}

\date{January 2013 \\ 
\vspace{6pt} 
\normalsize{Final version to be published in \\
{\it Stochastic Processes and their Applications}}}

\maketitle

\begin{abstract}
  We consider random walks in random environments on $\Z^d$. Under a
  transitivity hypothesis that is much weaker than the customary
  ellipticity condition, and assuming an absolutely continuous
  invariant measure on the space of the environments, we prove the
  ergodicity of the annealed process w.r.t.\ the dynamics ``from the
  point of view of the particle''. This implies in particular that the
  environment viewed from the particle is ergodic. As an example of
  application of this result, we prove a general form of the quenched
  Invariance Principle for walks in doubly stochastic environments
  with zero local drift (martingale condition).

  \medskip\noindent
  MSC 2010: 60G50, 60K37, 37A50 (37A20, 60G42, 60F17). 
\end{abstract}

\section{Introduction}
\label{sec-intro}

In this note we investigate \rw s in \ra\ \en s (RWREs) on $\Z^d$,
i.e., $\Z^d$-valued Markov chains defined by the transition matrix
$p(\om) = ( p_{xy}(\om) )_{x,y \in \Z^d}$, where $\om$ is a random
parameter ranging in a complete \pr\ space $(\Omega, \Pi)$. (In the
remainder we will refer to either $\om$ or $p(\om)$ as the
\emph{\en}.)

Although the precise nature of $\Omega$ is irrelevant, a natural
choice is $\Omega = (\mathcal{S}_{K,\gamma})^{\otimes \Z^d}$, where,
for $K,\gamma > 0$, $\mathcal{S}_{K,\gamma}$ is the space of all \pr\
distributions $\om_o = ( \om_{oy} )_{y \in \Z^d}$ on $\Z^d$, such that
$\om_{oy} \le K |y|^{-d-\gamma}$. By tightness \cite{b2},
$\mathcal{S}_{K,\gamma}$ is compact in the weak-* topology, which can
be metrized, e.g., by the total variation distance between two
distributions. Hence, by Tychonoff's Theorem and a standard argument,
$\Omega$ is also compact and metrizable. (This is important in case
one needs to construct suitable invariant \me s.)  For an element $\om
= ( \om_x )_{x \in \Z^d} \in \Omega$, where $\om_x = ( \om_{xy} )_{y
  \in \Z^d} \in \mathcal{S}_{K,\gamma}$, $p(\om)$ is defined by
$p_{xy}(\om) := \om_{x,y-x}$.

$\Omega$ is acted upon by $\Z^d$ via the group of $\Pi$-automorphisms
$( \tau_z )_{z \in \Z^d}$, such that
\begin{equation}
  \label{zd-action}
  p_{xy}(\tau_z \om) = p_{x+z,y+z}(\om).
\end{equation}
(In the representation above, $\tau_z$ is defined as $(\tau_z \om)_x
:= \om_{x+z}$.) Because of this, there is no loss of generality in
requiring that the walk always starts at $0$. The \rw\ (RW) in the
\en\ $p(\om)$ is then be defined as the Markov chain $( X_n )_{n \in
  \N}$ on $\Z^d$ whose law $P_\om$ is uniquely determined by
\begin{eqnarray}
  && P_\om (X_0 = 0) = 1; \label{rw-d1} \\
  && P_\om (X_{n+1} = y \,|\, X_n = x) = p_{xy}(\om). \label{rw-d2}
\end{eqnarray}
The complete \ra ness of the problem is accounted for by the
\emph{annealed} (or \emph{averaged}) law, which is defined on
$(\Z^d)^\N \times \Omega$ via
\begin{equation}
  \label{bbp}
  \mathbb{P} (E \times B) := \int_B \Pi(d\om) \, P_\om(E),
\end{equation}
where $E$ is a Borel set of $(\Z^d)^\N$ (the latter being the space of
the \tr ies, where $P_\om$ is defined) and $B$ is a measurable set of
$\Omega$. A natural dynamics that can be defined on this process is
the one induced by the map $\mathcal{F}: (\Z^d)^\N \times \Omega
\longrightarrow (\Z^d)^\N \times \Omega$, given by
\begin{equation}
  \label{calf}
  \mathcal{F} \left( (x_n)_{n\in\N} \,, \om \right) := \left( 
  (x_{n+1})_{n\in\N} \,, \tau_{x_1} \om \right).
\end{equation}
As is apparent, the first component of $\mathcal{F}$ updates the \tr y
of the RW to the next time, while the second component updates the
\en\ as seen by the \ra\ walker, or particle. This dynamics may thus
be called `the point of view of the particle' (PVP) for the annealed
process.

Another process of great importance, which is directly related to the
above, is the so-called `\en\ viewed from the particle' (EVP). It can
be defined independently as the Markov chain $(\Omega_n)_{n\in\N}$ on
$\Omega$, with law $\mathcal{P}_\Pi$, such that:
\begin{eqnarray}
  && \mathcal{P}_\Pi (\Omega_0 \in B) = \Pi (B); \label{evp-d1} \\
  && \mathcal{P}_\Pi (\Omega_{n+1} = \om \,|\, \Omega_n 
  = \omega') = \sum_{y: \, \tau_y \om' = \om} p_{0y}(\om'). \label{evp-d2}
\end{eqnarray}
(The annoying notation whereby $\Omega_n$ denotes an element of
$\Omega$ will be used only once more in this note, in the statement
of (A1).)

\bigskip\noindent 
\textsc{Notational convention.}\ Throughout the paper, the dependence
of whichever quantity (e.g., $p$) on $\om$ will not be explicitly
indicated when there is no risk of confusion.  
\medskip

When studying the stochastic properties of a RWRE, say, proving the
quenched Central Limit Theorem (CLT) (i.e., the CLT relative to
$P_\om$, for $\Pi$-a.a.\ $\om \in \Omega$), the state-of-the-art
technique requires one to face more or less three hurdles:
\begin{enumerate}
\item showing the existence of a steady state for the EVP that is
  absolutely continuous w.r.t.\ the original \me;
\item proving certain basic \erg\ properties of the steady state;
\item controlling the error term between the random walk and an
  approximating martingale.
\end{enumerate}

An assumption that is almost always made is the \emph{ellipticity} of
the \en. We state two rather general versions it may come in.

\begin{definition}
  \label{def-ell}
  A \ra\ \en\ $(\Omega, \Pi)$ is called \emph{elliptic} if there
  exists a basis $\Lambda$ of $\Z^d$ such that, for $\Pi$-a.a.\ $\om
  \in \Omega$ and all $e \in \Lambda$, $p_{0e} > 0$. It is called
  \emph{uniformly elliptic} if there also exists $\eps > 0$ such that,
  in the same cases as before, $p_{0e} \ge \eps$.
\end{definition}
To the author's knowledge, at least within the scope of non-ballistic
RWREs, only a few results \cite{be, ss, bb, mp, bed} do not assume
ellipticity. In general, uniform ellipticity is required---although
recent work focuses on non-uniformly elliptic \sy s; cf.\ \cite{m, gz,
sa} and references therein.

There are reasons to consider the ellipticity condition rather
unsatisfactory; for example, uniform ellipticity is a deterministic
condition in a probabilistic problem. The purpose of this note is to
show that, among the three hurdles mentioned above, ellipticity might
only be important for the first one. Certainly it is not needed for
the second. (The third hurdle is outside the scope of this article,
but the indication is that ellipticity is not particularly useful
there either.  See also Section \ref{sec-bis}.)

We replace ellipticity with the hypothesis that a relevant set of \en
s are \emph{transitive}, i.e., a walker starting at the origin goes
anywhere with positive probability. We call this \emph{partial
transitivity}.  (In actuality, our hypotheses are more general than
that; cf.\ (A3)-(A5) below.) We say that \emph{(total) transitivity}
holds if a.e.\ \en\ is transitive.  A certain form of transitivity,
together with the \erg ity of the \en\ w.r.t.\ translations, seems
like a natural, somewhat minimal, assumption, if we want the particle
to ``test'' every \en. A discussion on the merits of our hypothesis is
given below.

In a general setting, we show that, if an absolutely continuous steady
state is given on the space of the \en s (that is, if hurdle 1 is
taken care of), and if its support contains transitive \en s, the
steady state and the original state are equivalent \me s, and the PVP
process, relative to either, is \erg\ (Theorem \ref{thm-main}). The
same is true for the EVP process (Corollary \ref{cor-evp}).

The simplest \en s for which the above hypothesis holds are the doubly
stochastic \en s, for which $\Pi$ is automatically invariant in the
right sense. If we further assume the RW to be a martingale (doing
away with hurdle 3 as well) one gets the Invariance Principle (IP)
too (Theorem \ref{thm-qip}).

\bigskip\noindent
\textsc{Assumptions.}\ The RWRE satisfies the following:

\begin{itemize}
\item[(A1)] \emph{No deterministic walks}. Denote by $\delta_{xy}$ the
  Kronecker delta in $\Z^d$ and set
  \begin{displaymath}
    \Omd := \rset{\om \in \Omega} {\exists y_o \in \Z^d 
    \mbox{ such that } p_{0y} = \delta_{y_o y}, \forall y \in \Z^d}.
  \end{displaymath}
  (This is the set of \en s that have a deterministic jump at the
  origin.)  We have that, $\mathbb{P}$-almost surely,
  $(\Omega_n)_{n\in\N} \not\subset \Omd$. In other words, in a.e.\
  \en, the \rw\ starting at the origin is not deterministic.

\item[(A2)] \emph{Decaying transition probabilities}. There exist $K,
  \gamma>0$ such that, almost surely,
  \begin{displaymath}
    p_{xy} \le K |y-x|^{-d-\gamma}.
  \end{displaymath}

\item[(A3)] \emph{Ergodicity}. There is a subgroup $\Gamma \subseteq
  \Z^d$ such that $(\Omega, \Pi, ( \tau_z )_{z \in \Gamma} )$ is
  \erg. (In particular, the \ra\ \en\ is \erg\ w.r.t.\ the whole group
  of translations.)

\item[(A4)] \emph{Partial transitivity}. Let $\Gamma$ be the same as
  in (A3).  Define
  \begin{displaymath}
    \Omt := \rset{\om \in \Omega} { \forall y \in
    \Gamma, \ \exists n = n(\om,y) \mbox{ such that } p^{(n)}_{0y} > 0 },
  \end{displaymath}
  where $p^{(n)}_{0y} := \sum_{x_1, \ldots, x_{n-1}} p_{0x_1} \,
  p_{x_1 x_2} \cdots\, p_{x_{n-1} y}$. (This is the set of transitive
  \en s.) It must be that $\Pi(\Omt) > 0$.

\item[(A5)] \emph{Absolute continuity of steady state}. There exists a
  \pr\ \me\ $\ime$ on $\Omega$, absolutely continuous w.r.t.\ $\Pi$,
  such that the EVP process $\mathcal{P}_{\ime}$, defined as in
  (\ref{evp-d1})-(\ref{evp-d2}), is stationary. Furthermore, $\Pi(
  \mathrm{supp}_\Pi \, \ime \cap \Omt ) > 0$; equivalently, $\ime(
  \Omt ) > 0$.
\end{itemize}
Notice that there is a trade-off between (A3) and (A4): if $\Gamma$
gets smaller, thus making (A3) stronger, then (A4) becomes easier to
verify, and viceversa. In particular, for an i.i.d.\ \en\ (namely, the
stochastic vectors $p_x = ( p_{xy} )_{y \in \Z^d}$ are i.i.d.\ in
$x$), one only need verify partial transitivity w.r.t.\ $\Gamma = \Z
y_o$, for some $y_o \in \Z^d$. In any event, $\Gamma = \Z^d$ is a
reasonable choice for many applications. For further observations on
(A4), see Remark \ref{rk-pt-trans}.

\begin{remark}
  In (A5), the hypothesis of absolute continuity could be replaced by
  the (slightly) weaker condition that the steady state $\ime$ is
  non-singular w.r.t.\ $\Pi$. The two conditions are actually
  equivalent in our case. In fact, it can be seen (cf.\ Section
  \ref{sec-dsy}) that the evolution of $\Pi$ in the EVP process is
  absolutely continuous w.r.t.\ $\Pi$, which implies that, if $\ime$
  decomposes into an absolutely continuous \me\ and a singular \me,
  both of them are invariant for the process.
\end{remark}

\noindent
\textsc{Results.}\ These are our main results:

\begin{theorem}
  \label{thm-main}
  Under assumptions \emph{(A1)-(A5)}, 
  \begin{itemize}
  \item[(a)] the \me s $\Pi$ and $\ime$ are equivalent (i.e., mutually 
    absolutely continuous);

  \item[(b)] if $\ima$ is the annealed law relative to $\ime$ (i.e.,
    the \me\ defined by \emph{(\ref{bbp})} with $\ime$ in lieu of
    $\Pi$), then $\ima$ is stationary and \erg\ for the dynamics
    induced by $\mathcal{F}$ on the annealed process;
    
  \item[(c)] $\Pi(\Omt) = 1$.  
  \end{itemize}
\end{theorem}

\begin{corollary}
  \label{cor-evp}
  The EVP with initial state $\ime$ is \erg.
\end{corollary}
Another easy corollary of Theorem \ref{thm-main} concerns the
ballisticity of the RW. In order to state it, we introduce the
\emph{mean displacement} (or \emph{local drift}) at the origin, for
the \en\ $\om$. This is the \fn\ $\ld: \Omega \longrightarrow \Z^d$
given by
\begin{equation}
  \label{l-dr}
  \ld(\om) := \sum_{y \in \Z^d} p_{0y}(\om) \, y,
\end{equation}
which is surely well-defined if $\gamma$ in (A2) is bigger than 1.
\begin{corollary}
  \label{cor-ball}
  If the local drift (\ref{l-dr}) is well-defined in $L^1 (\Omega,
  \Pi)$, then, for $\Pi$-a.e.~\en\ $\om$,
  \begin{displaymath}
    \lim_{n \to \infty} \frac{X_n}n = \int_\Omega \ime(d\om') \,
    \ld (\om'),
  \end{displaymath}
   $P_\om$-almost surely.
\end{corollary}

\bigskip\noindent 
\textsc{Discussion.}\ Theorem \ref{thm-main} is a ``soft result'', in
the sense that it is very general and not very deep. Certainly it is
not surprising, as analogous statements were already known in several
specific cases, e.g., \cite{la, la2, bk, ko, ss, bb, mp, m, bad, gz,
  sa, bed}.  But this is actually the point of this note: to show that
the \erg ity of the PVP is a general result that need not be proved
every time, provided one has an absolutely continuous invariant
measure for the EVP, and transitivity somewhere in its support.  In
other words, hurdle 2 is not a hurdle.

It might be worthwhile to point out that Theorem \ref{thm-main} is not
obvious from Koslov's 1985 paper \cite{k} (at least not to this
author).

Recently, Berger and Deuschel \cite{bed} proved the quenched IP for
RWs in i.i.d., nearest-neighbor, balanced \en s (in the sense of
\cite{la, la2}), under the assumption of \emph{genuine
$d$-dimensionality} (\gdd). The latter means that, for every $e \in
\Z^d$ with $|e|=1$, $\Pi( p_{0e}(\om) > 0 ) >0$. Clearly, this is
weaker than our transitivity; strictly weaker, in fact, as we discuss
below.

In light of this result, one might think that \gdd\ is a better
condition than transitivity and would produce a more general version
of Theorem \ref{thm-main}. We claim that this is not really the
case. First we argue that, within the scope of \cite{bed}, \gdd\ is
similar in spirit to transitivity. Then, by way of a counterexample,
we show that Theorem \ref{thm-main} could not hold in its generality
if transitivity were replaced by \gdd

Under the assumptions of \cite{bed}, it is easy to find examples
where, almost surely, there are sites that the particle cannot visit
(see, e.g., Fig.\ 3 of \cite{bed}). It turns out, however, that a.e.\
\en\ has a \emph{sink}, that is, a subset of $\Z^d$ that each \tr y
enters, in finite time, and never leaves. The sink is unique,
unbounded, and transitive. The \en s seen from the sink make up the
support of the steady state $\ime$ (whose $\Pi$-\me\ may be less than
1) and its \erg\ properties depend on the fact that it is transitive.

As for the other point---disproving Theorem \ref{thm-main} under the
hypothesis of \gdd---we have just shown that assertion \emph{(a)} may
not be true. This is no big trouble, as long as one proves that, after
a controlled time, the EVP falls in the support of $\ime$ (this is the
case in \cite{bed}). But there are examples in which assertion
\emph{(b)} fails as well.

For instance, let $( \xi_j )_{j \in \Z}$ be i.i.d.\ \ra\ variables
taking the values (labels) A or B, both with positive
probability. Each realization of this Bernoulli chain defines an \en\
on $\Z^2$ by assigning to the site $x = (j,k)$ the same label as the
variable $\xi_j$. Thus, the \ra\ \en\ is constant for the vertical
translations and \erg\ for the horizontal ones, whose subgroup we
denote $\Gamma$; cf.\ (A3). It is then \erg. Prescribe that, when the
particle visits a site of type A, it has \pr\ 1/2 to make the next
move to the left and \pr\ 1/2 to make it to the right; when it visits
a site of type B, the rule is analogous, but for up/down
moves. Clearly, \gdd\ holds but transitivity (almost) never does, even
for $\Gamma$; cf.\ (A4).  The associated PVP and EVP processes have
quite trivial dynamics: in particular, the elementary steady states
for the EVP are atoms, corresponding to \en s whose entire vertical
axis is labeled B. This shows that an absolutely continuous steady
state $\ime$ cannot make those processes \erg.

\begin{remark}
  The fact that the above \en\ is balanced with nearest-neighbor
  jumps, and \gdd\ holds, proves that a strong mixing property is
  essential for the result of \cite{bed}.
\end{remark}

Our proofs are based on a convenient representation of the RWRE as a
\pr-preserving \dsy\ which, roughly speaking, is a ``\me d family'' of
one-dimensional Markov maps. Each map embodies the dynamics of one
\ra\ jump, and thus contains only local information. We will see that
this \dsy\ is isomorphic to the annealed process. In any case, Section
\ref{sec-dsy} should convince the reader that it is natural to call
this object `the \dsy\ for the point of view of the particle'; in
short, \emph{PVP \dsy}.

\skippar

The exposition is organized as follows: In Section \ref{sec-dsy} we
introduce the \dsy\ and find a suitable invariant measure for it. In
Section \ref{sec-erg} we prove its ergodicity, which is equivalent to
Theorem \ref{thm-main} and implies its corollaries. In Section
\ref{sec-bis} we consider the example of the doubly stochastic RWs,
proving an improved version of the quenched IP for doubly stochastic
martingales.

\bigskip\noindent 
\textsc{Acknowledgments.}\ I am grateful to Alessandra Bianchi, Firas
Rassoul-Agha, Frank den Hollander, Luc Rey-Bellet, Vladas Sidoravicius
and Stefano Olla for useful discussions. I also thank an anonymous
referee for a careful reading of the early versions of the manuscript.
This work was partially supported by the FIRB-``Futuro in Ricerca''
Project RBFR08UH60 (MIUR, Italy).

\section{The PVP dynamical system}
\label{sec-dsy}

Let us fix an enumeration $( d_i )_{i \in \Z_+}$ of $\Z^d$. For $\om
\in \Omega$ and $i \in \Z_+$, we define
\begin{equation}
  \label{qi}
  q_i(\om) := p_{0d_i}(\om) \\ 
\end{equation}
Certainly, $\sum_i q_i(\om) = 1$. We then set $a_0(\om) := 0$ and,
recursively on $i$,
\begin{eqnarray}
  \label{ai}
  && a_i(\om) := a_{i-1}(\om) + q_i(\om) \\
  \label{ii}
  && I_i(\om) := [a_{i-1}(\om) \,,\, a_i(\om)).
\end{eqnarray}
Clearly, $\{ I_i(\om) \}_{i \in \Z_+}$ is a partition of $I := [0,1)$.
For $(s,\om) \in I \times \Omega$, denote by $i(s,\om)$ the unique $i$
such that $s \in I_i(\om)$. We define the \fn\ $\phi: I \times \Omega
\longrightarrow I$ via
\begin{equation}
  \label{def-phi}
  \phi(s,\om) := \frac{ s - a_{i(s,\om) - 1} (\om) }
  { q_{i(s,\om)} (\om) }.
\end{equation}
(The definition above is well-posed because, if $i$ is such that
$q_i(\om) = 0$, there is no $s$ such that $i(s,\om) = i$.)  By
construction, $\phi(\cdot, \om)$ is the perfect Markov map $I
\longrightarrow I$ relative to the partition $\{ I_i(\om) \}$, namely,
each branch of this map is affine and its image is $I$.
Finally, we denote $D(s,\om) := d_{i(s,\om)}$.

The main technical tool of this paper is the map $\ma: \ps
\longrightarrow \ps$, defined on $\ps := I \times \Omega$ by
\begin{equation}
  \label{def-ma}
  \ma (s,\om) := \left( \phi(s,\om) , \tau_{D(s,\om)} \om \right).
\end{equation}
We endow $\ps$ with either the \pr\ \me\ $\mu := m \otimes \Pi$ or $\im
:= m \otimes \ime$, where $m$ is the Lebesgue \me\ on $I$.

What this \dsy\ has to do with our RWRE is presently explained.  Let
us recall the notational convention whereby the dependence on $\om$ is
not always indicated. Fix $\om \in \Omega$ and consider a \ra\ $s \in
I$ w.r.t.\ $m$. We have that $D(s,\om) = d_i$ \iff $s \in I_i$, and
this occurs with \pr\ $m(I_i) = q_i$. In terms of our RW, this is
exactly the \pr\ that a particle placed in the origin of $\Z^d$,
endowed with the \en\ $p(\om)$, jumps by a quantity $d_i$. Then, back
to the \dsy, condition the \me\ $m$ to $I_i$. Calling $(s_1,\om_1) :=
\ma (s,\om)$, we see that, upon conditioning, $s_1$ ranges in $I$ with
law $m$. Therefore, in a sense, the variable $s$ (which we may call
the \emph{internal variable}) has ``refreshed'' itself. Furthermore,
$\om_1 = \tau_{D(s,\om)} \om = \tau_{d_i} \om$ is the translation of
$\om$ in the opposite direction to $d_i$; cf.\ (\ref{zd-action}).
Hence we can imagine that we have reset the \sy\ to a new initial
condition $(s_1,\om_1)$, corresponding to the particle sitting in $0
\in \Z^d$ and subject to the \en\ $p(\om_1)$. Applying the same
reasoning to $(s_2,\om_2) := \ma (s_1,\om_1)$, and so on, shows that
we are following the motion ``from the point of view of the
particle''. We thus call the above the `PVP \dsy'.

In any case, it should be clear that the stochastic process $( X_n
)_{n \in \N}$, with $X_0 \equiv 0$ and, for $n \ge 1$,
\begin{equation}
  \label{def-xn}
  X_n(s,\om) := \sum_{k=0}^{n-1} D \circ \ma^k(s,\om).
\end{equation}
is precisely the RW in the \en\ $p(\om)$, provided that $\om$ is
regarded as a fixed parameter. To emphasize this point, we
occasionally write $X_{n,\om}(s) := X_n(s,\om)$. $( X_{n,\om} )$ is
called the `quenched \tr y', and it is a Markov chain. If both $s$ and
$\om$ are considered random, w.r.t.\ $\mu$, then (\ref{def-xn})
defines the `annealed \tr y'. This is not a Markov chain and, by the
definition of $\mu$ and (\ref{bbp}), it is none other than the RWRE of
Section \ref{sec-intro} with law $\mathbb{P}$. For a formal relation
between the annealed process and the PVP \dsy\ see Proposition
\ref{prop-iso}.

\begin{proposition}
  \label{prop-im}
  The \me\ $\im$ is preserved by $\ma$.
\end{proposition}
We need the following lemma:
\begin{lemma}
  \label{lem-ss}
  For every measurable set $B \subseteq \Omega$,
  \begin{displaymath}
    \ime(B) = \sum_{i \in \Z_+} \int_{\tau_{-d_i(B)}} \ime(d\om') \, 
    q_i(\om').
  \end{displaymath}
\end{lemma}

\proofof{Lemma \ref{lem-ss}} Via (\ref{qi}) and (A2), we observe that
the series
\begin{equation}
  \label{steady5}
  \sum_{i \in \Z_+} q_i(\om') = \sum_{x \in \Z^d} p_{0x}(\om'),
\end{equation}
is uniformly bounded in $\om'$. Thus, as we will do more than once
presently, it is correct to interchange the above summation with an
integration over $\Omega$, if it is relative to a \pr\ \me.

Using again (\ref{qi}), the transition kernel of the EVP process
(\ref{evp-d2}) can be written as
\begin{equation}
  \label{steady10}
  \mathcal{K} ( \om', \, \cdot \,) = \sum_{i \in \Z_+} q_i(\om') \, 
  \delta_{\tau_{d_i} \om'},
\end{equation}
where the Dirac delta on the r.h.s.\ is thought of as a \me. In other
words, for a measurable set $B$,
\begin{equation}
  \label{steady20}
  \mathcal{K} (\om', B) = \int_\Omega \mathcal{K} (\om', d\om) \, 
  1_B(\om) = \sum_{i \in \Z_+} q_i(\om') \, 1_{\tau_{-d_i(B)}} (\om').
\end{equation}
Thus, the hypothesis on $\ime$ from (A5), namely,
\begin{equation}
  \label{steady30}
  \ime(B) = \int_\Omega \ime(d\om') \, \mathcal{K} (\om', B),
\end{equation}
reads precisely as in the statement of the lemma.  
\qed

\proofof{Proposition \ref{prop-im}} It is enough to prove that
$\im(\ma^{-1} A) = \im(A)$ for all sets of the type $A = [b,c) \times
B$, where $B$ is a measurable set of $\Omega$.

By direct inspection of the map (\ref{def-ma}), we can write $\ma^{-1}
A = \bigcup_{i \in \Z_+} A'_i$, where
\begin{equation}
  \label{meas10}
  A'_i := \rset{ (s',\om') } {\om' \in \tau_{-d_i}(B), \, s' \in 
  [ a_i(\om') + q_i(\om') b \,,\, a_i(\om') + q_i(\om') c) };
\end{equation}
cf.\ (\ref{ai}), (\ref{def-phi}). Thus,
\begin{equation}
  \label{meas20}
  \im(A'_i) = \int_{\tau_{-d_i}(B)} \hspace{-4pt} \ime(d\om') \,
  q_i(\om') (c-b).
\end{equation}
The sets $A'_i$ are pairwise disjoint because, by construction, they
belong to different level sets of the \fn\ $D$.  Therefore, by Lemma
\ref{lem-ss},
\begin{equation}
  \label{meas30}
  \im(\ma^{-1} A) = (c-b) \sum_{i \in \Z_+} \int_{\tau_{-d_i}(B)} 
  \hspace{-4pt} \ime(d\om') \, q_i(\om') = (c-b) \, \ime(B) = \im(A),
\end{equation}
which is what we wanted to prove.
\qed

Let us introduce a convenient notation that will be used throughout
the paper: For $(s,\om) \in \ps$ and $k \in \N$, denote
\begin{equation}
  \label{sk-omk}
  (s_k,\om_k) := \ma^k(s,\om).
\end{equation}

\begin{proposition}
  \label{prop-ac}
  The \me s $\im$ and $\mu$ are equivalent or, which is the same, the
  \me s $\ime$ and $\Pi$ are equivalent.
\end{proposition}

\proof Denote $\sime := \mathrm{supp}_\Pi \, \ime$,
$\Omega_\mathrm{*t} := \sime
\cap \Omt$, and
$\sime^c := \Omega \setminus \sime$. By (A5), $\Pi(\Omega_\mathrm{*t}) > 0$
and $\Pi (\sime^c) <
1$. Suppose, by way of contradiction, that $\Pi (\sime^c) > 0$ as
well.

By (A3), for $\Pi$-a.e.\ $\om \in \Omega_\mathrm{*t}$, there exists $y_\om \in
\Gamma$ such that
\begin{equation}
  \label{im30}
  \tau_{y_\om} \om \in \sime^c.
\end{equation}
To each such $\om$ (excluding at most a $\Pi$-null set) we apply (A4)
and its interpretation in terms of the PVP \dsy: there exists a
positive integer $n = n(\om, y_\om)$ such that
\begin{equation}
  \label{im35}
  J_\om := \rset{s \in I} {X_n (s,\om) = y_\om}
\end{equation}
has \me\ $p_{0 y_\om}^{(n)} > 0$. By definition, $\forall s \in
J_\om$,
\begin{eqnarray}
  \ma^n (s,\om) &=& (s_n, \tau_{D(s_{n-1},\om_{n-1})} \circ \cdots 
  \circ \tau_{D(s,\om)} \, \om) \nonumber \\
  &=& (s_n, \tau_{X_n(s,\om)} \, \om) \nonumber \\
  \label{im40}
  &=& (s_n, \tau_{y_\om} \om) \in I \times \sime^c,
\end{eqnarray}
having used (\ref{sk-omk}), (\ref{def-xn}), (\ref{im35}) and finally
(\ref{im30}). If we define
\begin{equation}
  \label{im45}
  A := \bigcup_{\om \in \Omega_\mathrm{*t}} J_\om \times \{ \om \},
\end{equation}
we have $\mu(A) > 0$ and, via (\ref{im40}),
\begin{equation}
  \label{im50}
  A \subset \bigcup_{n \ge 1} T^{-n} (I \times \sime^c).
\end{equation}
The definition of $\sime^c$ implies that $\im (I \times \sime^c) = 0$,
hence, since $\im$ is $\ma$-invariant, $\im(A) = 0$. Finally, since
$\mu$ and $\im$ are equivalent on $I \times \sime \supset A$,
$\mu(A)=0$, which contradicts a previous statement.  
\qed

For a graphic representation of the above proof, see
Fig.~\ref{rwwe-fig} and its caption.

\fig{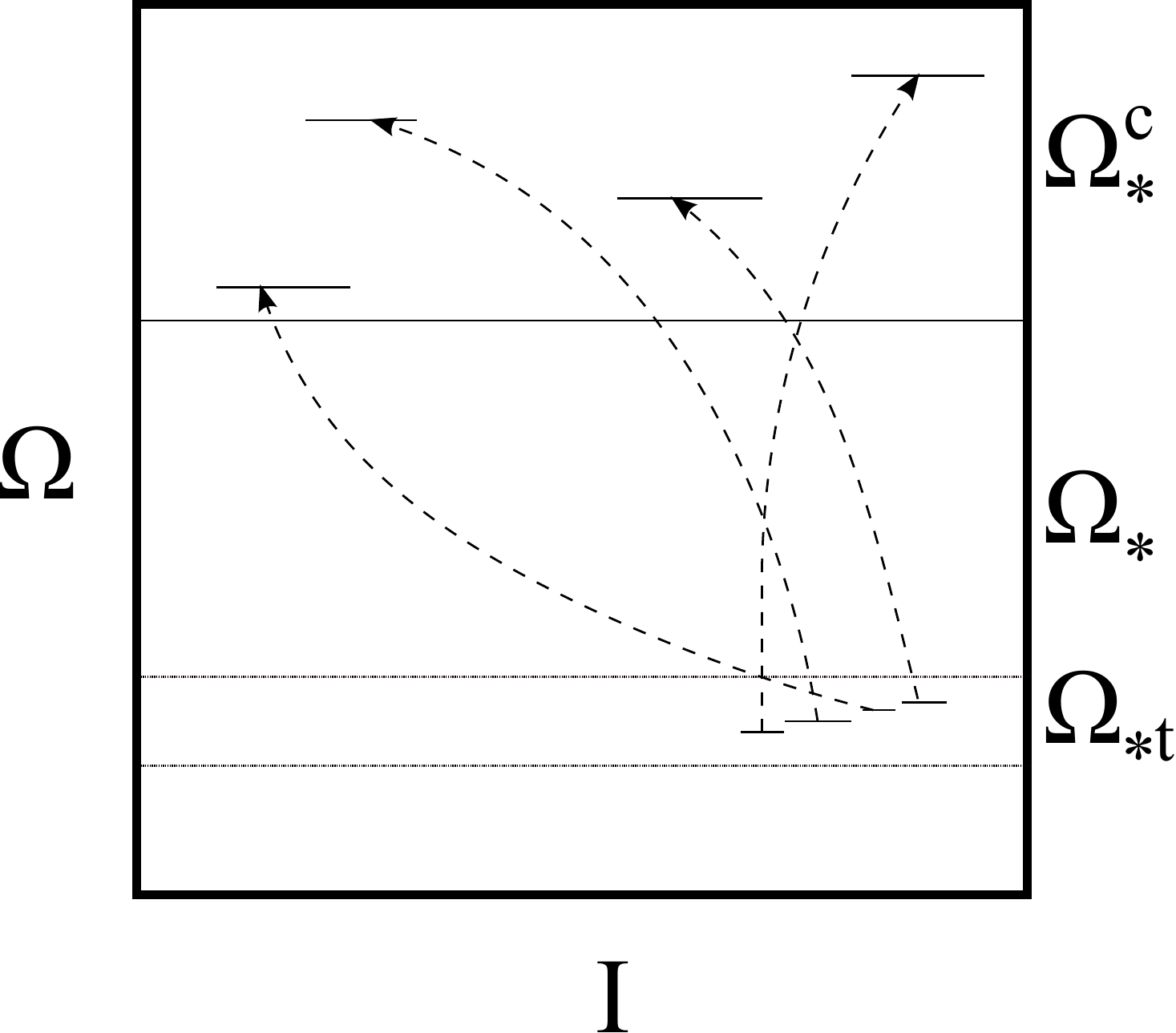}{7cm} {Representation of $\ps = I \times \Omega$ as a
  square for the proof of Proposition \ref{prop-ac}. By the \erg ity
  of the \ra\ \en\ and partial transitivity, in a.e.\ horizontal
  segment (later denoted by \emph{horizontal fiber}) of
  $\Omega_\mathrm{*t} \subseteq \Omt$, one can find a subsegment
  (denoted by $J_\om \times \{ \om \}$ in the proof) that the dynamics
  takes ``everywhere in $\ps$'' that has been established. In this
  way, one carries a positive $\mu$-\me\ into $I \times \sime^c$, the
  complement of supp$\,\im$, which therefore cannot be $\mu$-null.}

\section{Ergodicity}
\label{sec-erg}

In this section we will prove the \erg ity of $(\ps, \im, \ma)$.  For
this, we need to introduce some more notation and establish a few
lemmas.

Given a positive integer $n$ and a multi-index $\mathbf{i} := (i_0,
i_1, ..., i_{n-1}) \in \Z_+^n$, we set
\begin{equation}
  \label{iib}
  I_\mathbf{i} (\om) := \rset{s \in I} {D \circ \ma^k (s,\om) =
  d_{i_k},\: \forall k=0, \ldots, n-1}.
\end{equation}
For $n=1$ this reduces to definition (\ref{ii}). It is easy to
ascertain that $\{ I_\mathbf{i} \}_{\mathbf{i} \in \Z_+^n}$ is a
partition of $I$ into countably many (possibly empty) right-open
intervals, each of which corresponds to one of the realizations of the
RW $( X_{k,\om} )_{k=0}^n$ (relative to the \en\ $\om$) in such a way
that $m (I_\mathbf{i})$ is the \pr\ of the corresponding
realization. In analogy with the previous notation, we denote by
$\mathbf{i}_n(s,\om)$ the index of the element of the partition which
contains $s$.

Furthermore, let us call \emph{horizontal fiber} of $\ps$ any segment
of the type $I_\om := I \times \{ \om \}$, and indicate with $m_\om$
the Lebesgue \me\ on it. Finally, we denote by $I_{\om,\mathbf{i}}$ the
subset of $I_\om$ corresponding to $I_\mathbf{i} (\om)$ via the
natural isomorphism $I \cong I_\om$.

\begin{lemma}
  \label{lem-hyp}
  For a.a.\ $(s,\om) \in \ps$, $m( I_{\mathbf{i}_n(s,\om)} )$ vanishes
  exponentially fast, as $n \to \infty$.
\end{lemma}

\proof Set
\begin{equation}
  \label{erg5}
  f(s,\om) := \log q_{i(s,\om)}^{-1} (\om) = -\log m(I_{i(s,\om)}
  (\om)).
\end{equation}
Clearly $f(s,\om) \ge 0$, with $f(s,\om) = 0$ \iff $\{ I_i(\om) \}$ is
the trivial partition of $I$, \iff $\om \in \Omd$; cf.\ (A1).  The
Birkhoff average of $f$,
\begin{equation}
  \label{erg10}
  f^+ (s,\om) := \lim_{n \to \infty} \, \frac1n \, \sum_{k=0}^{n-1} \,
  f (s_k,\om_k) ,
\end{equation}
is non-negative as well. Set $A := \rset{(s,\om)} {f^+(s,\om) =
0}$. As a level set of an invariant \fn, $A$ is invariant mod
$\im$. We claim that $\im(A)=0$. If not, we can apply one of the
assertions of Birkhoff's Theorem to the \me-preserving \dsy\ $(A, \im,
\ma_{|A})$ and conclude that $\int_A d\im \, f = \int_A d\im \, f^+ =
0$. Therefore, $f(s,\om) = 0$, for a.a.\ $(s,\om) \in A$. In other
words, $A \subseteq I \times \Omd$ mod $\im$. Since $A$ is invariant,
the \o\ of a.e.\ point in $A$ is contained in $I \times \Omd$,
implying that $(\om_n)_{n \in \N} \subset \Omd$ for a positive \me\ of
initial conditions---in contradiction with (A1).

So $f^+>0$ almost everywhere. On the other hand, from earlier
considerations, it is easy to verify that, for $n \ge 1$,
\begin{equation}
  \label{erg20}
  m \left( I_{\mathbf{i}_n(s,\om)} (\om) \right) = \prod_{k=0}^{n-1} 
  q_{i(s_k,\om_k)} (\om_k) = \exp \left( - \sum_{k=0}^{n-1} 
  f (s_k,\om_k) \right).
\end{equation}
Due to the almost sure positivity of (\ref{erg10}), the exponent in
the rightmost term above is asymptotically linear in $n$, for a.a.\
$(s,\om)$, which yields the assertion.  
\qed

\begin{proposition}
  \label{prop-iso}
  As \dsy s on \pr\ spaces, $(\ps, \mu, \ma)$ is isomorphic to $(
  (\Z^d)^\N \times \Omega, \mathbb{P}, \mathcal{F})$, and $(\ps, \im,
  \ma)$ is isomorphic to $( (\Z^d)^\N \times \Omega, \ima,
  \mathcal{F})$.
\end{proposition}

\proof We hope the reader was already convinced in Section
\ref{sec-dsy} that the PVP \dsy\ describes exactly the annealed
process with the PVP dynamics.  On the other hand, Lemma \ref{lem-hyp}
provides the ingredients for a formal proof, which we just sketch
here.

For both pairs of \sy s, a natural isomorphism $\Phi: \ps
\longrightarrow (\Z^d)^\N \times \Omega$ is given by
\begin{equation}
  \label{phi}
  \Phi (s,\om) := \left( (X_n (s,\om))_{n \in \N} \,, \om \right);
\end{equation}
cf.\ (\ref{def-xn}). One sees that $\Phi$ is almost-everywhere
bijective because of the following: By Lemma \ref{lem-hyp}, a.e.\
$(s,\om) \in \ps$ is the unique intersection point of the nested
sequence of right-open intervals $( I_{\om, \mathbf{i}_n(s,\om)}
)_{n\in\N}$, i.e., is uniquely determined by the sequence
$(\mathbf{i}_n(s,\om))$, equivalently, by the realization $(X_n
(s,\om))$ of the RW. Viceversa, for an \en\ $\om$, every realization
of the walk determines a nested sequence of intervals which, except
for a null set of realizations, gives a point $(s,\om) \in
I_\om$. (Since the intervals are right-open, this correspondence is
ill-defined at the endpoints of all such intervals. But this amounts
to a null set of points in $I_\om$ and a null set of realizations in
$(\Z^d)^\N$.)

Finally, it is clear by the considerations of Section \ref{sec-dsy}
that $\mu = \mathbb{P} \circ \Phi$, $\im = \ima \circ \Phi$ and $\ma =
\Phi^{-1} \circ \mathcal{F} \circ \Phi$.
\qed

\begin{lemma}
  \label{lem-erg}
  The \erg\ components of $(\ps, \im, \ma)$ contain whole horizontal
  fibers, that is, every invariant set is of the form $I \times B$,
  mod $\im$ (equivalently, mod $\mu$), where $B$ is a measurable
  subset of $\Omega$.
\end{lemma}

\proof The idea of the proof is that, since $\ma$ maps subintervals of
horizontal fibers onto whole horizontal fibers (the ``perfect Markov''
property of Section \ref{sec-dsy}), its repeated application stretches
any set horizontally to the full length of the horizontal fibers. An
invariant set must thus contain full fibers. Now for a formal proof.

Suppose the assertion of Lemma \ref{lem-erg} is false. There exists an
invariant set $A$ whose intersection with many horizontal fibers is
neither the full fiber nor empty, mod $m_\om$.  That is, for some
$\eps>0$, the $\ime$-\me\ of
\begin{equation}
  \label{erg30}
  B_\eps := \rset{\om \in \Omega} {m_\om (A \cap I_\om) \in [\eps,
  1-\eps]}
\end{equation}
is positive. By the Poincar\'e Recurrence Theorem and the Lebesgue
Density Theorem it is possible to pick $(s,\om) \in A \cap (I \times
B_\eps)$ that is recurrent to $I \times B_\eps$ and such that
$(s,\om)$ is a density point of $A \cap I_\om$, within $I_\om$,
relative to $m_\om$. We claim that there exist a sufficiently large
$n$ and a multi-index $\mathbf{i} \in \Z_+^n$ for which
\begin{equation}
  \label{erg35}
  m_\om( A \cap I_{\om,\mathbf{i}} ) > (1-\eps) \, m_\om(
  I_{\om,\mathbf{i}} )
\end{equation}
and 
\begin{equation}
  \label{erg40}
  \ma^n I_{\om,\mathbf{i}} = I_{\om_n} \subset I \times B_\eps. 
\end{equation}
In fact, among the infinitely many $n$ that verify $\ma^n(s,\om) \in I
\times B_\eps$, we can choose, by Lemma \ref{lem-hyp}, one for which
$I_{\om,\mathbf{i}_n (s,\om)}$ is so small that (\ref{erg35}) is
verified for $\mathbf{i} = \mathbf{i}_n (s,\om)$. The equality in
(\ref{erg40}) is true by the Markov property of $\phi(\cdot,
\om_{n-1}) \circ \cdots \circ \phi(\cdot, \om)$ (recall the notation
(\ref{sk-omk})).

It is no loss of generality to assume that the $n^\mathrm{th}$ iterate
of $m_\om$-a.e.\ point of $A \cap I_{\om,\mathbf{i}}$ remains in $A$
(in the choice of $(s,\om)$, use the invariance of $A$ mod $\im$ and
Fubini's Theorem).  Since the restriction of $\ma^n$ to
$I_{\om,\mathbf{i}}$ is linear, we deduce from
(\ref{erg35})-(\ref{erg40}) that $m_{\om_n} (A \cap I_{\om_n}) >
1-\eps$, which contradicts (\ref{erg30}), because $\om_n \in B_\eps$
by (\ref{erg40}).

Therefore, an invariant set mod $\im$ can only occur in the form $I
\times B$. The completeness of $\Pi$ (equivalently, $\ime$) implies
that $B$ is measurable, as in Lemma A.1 of \cite{l1} (cf.\ Lemma 3.4
of \cite{l3}).
\qed

\begin{remark}
  The techniques of Lemma \ref{lem-erg} (based on the fact that
  $\phi(\cdot, \om)$ is a piecewise-linear Markov map of the interval)
  easily imply that any $\ma$-invariant \me\ that is smooth along the
  horizontal fibers must be uniform on them, i.e., must be of the type
  $m \otimes \Pi'$.
\end{remark}

\begin{theorem}
  \label{thm-erg}
  $(\ps, \im, \ma)$ is \erg.
\end{theorem}

\proof Suppose the \sy\ is not \erg. By Lemma \ref{lem-erg}, we have
an invariant set $I \times B$, with $\Pi(B) \in (0,1)$. Either $B$ or
$B^c := \Omega \setminus B$ has a non-negligible intersection with
$\Omt$. Assume, without loss of generality, that the former is the
case.

The \pr\ $\mu' := \im (\,\cdot\, | I \times B)$ is $\ma$-invariant and
factorizes as $\mu' = m \otimes \Pi'$, where $\Pi' := \ime (\,\cdot\,
| B)$. Since $\Pi (\mathrm{supp}_\Pi \, \Pi' \cap \Omt ) > 0$, we can
apply Proposition \ref{prop-ac} with $\Pi'$ in the role of $\ime$. The
result contradicts the hypothesis $\Pi(B) \in (0,1)$ (cf.\
Fig.~\ref{rwwe-fig} with $B$ and $B^c$ in place of, respectively,
$\sime$ and $\sime^c$).
\qed

We can now easily prove the main result of the paper.

\medskip

\proofof{Theorem \ref{thm-main}} Assertion \emph{(a)} is Proposition
\ref{prop-ac}. Assertion \emph{(b)} follows from Proposition
\ref{prop-im}, Theorem \ref{thm-erg} and Proposition \ref{prop-iso}.

As for assertion \emph{(c)}, by Theorem \ref{thm-erg}, the \o\ of
$\im$-a.e.\ $(s,\om) \in \ps$ intersects $I \times \Omt$. This means
that, for $\ime$-a.e.\ $\om$, there exists $n_1$ such that, with
positive \pr\ (in the sense of $m_\om$), at time $n_1$ the particle
visits a site $x$ that is the origin of a transitive \en. From there
on, for any $y \in \Gamma$, there exists $n_2 = n_2(\tau_x \om, y)$
such that the particle has a positive \pr\ to be in $y$ after a time
$n_2$. By the Markov property of the \rw, this gives $p_{0y}^{(n_1 +
n_2)} (\om) > 0$.  In other words, using \emph{(a)}, $\Pi$-a.e.\
$\om$ is transitive.  
\qed

\begin{remark}
  \label{rk-pt-trans}
  Since assertion \emph{(c)} excludes---under (A1)-(A5)---partial but
  not total transitivity, one might opine that Theorem \ref{thm-main}
  would have been better formulated directly with the hypothesis of
  total transitivity. The current formulation, however, has some
  advantages.
  
  First off, it is stronger. And there are cases in which (A4) is easy
  to show, but it is not evident that $\Pi(\Omt) = 1$. For example,
  take an i.i.d.\ \en\ with a positive fraction of sites labeled C,
  where a C site $x$ has the property that $p_{xy} > 0$, $\forall y
  \in \Z^d$.  (A4) clearly holds, but it may not be clear that the
  other sites will make a.e.\ \en\ transitive (there might be sinks,
  or the like, see the discussion in the introduction).
  
  Second, Theorem \ref{thm-main} has the advantage that it can be
  easily adapted to the case in which the support of $\ime$ is
  strictly smaller than $\Omega$ (cf.\ \cite{bed}) and possibly there
  are more $\Pi$-absolutely continuous \erg\ states for the EVP.  Here
  is an example of such a result: For $B \subseteq \Omega$, say that
  the \en\ $\om$ is \emph{$B$-transitive} if, $\forall y \in \Gamma$
  with $\tau_{-y} \om \in B$, $\exists n=n(y,\om)$ such that
  $p_{0y}^{(n)} (\om) > 0$. In other words, $\om$ is $B$-transitive if
  the walker has a positive probability to go to any site (of
  $\Gamma$) where he sees an \en\ in $B$. (Thus, $\Omega$-transitive
  means transitive.)  Then, if (A1)-(A3) hold and $\ime$ as in (A5)
  exists, the non-null \erg\ components of $\ime$ are the maximal sets
  $B$ such that, for a fixed $B_0 \subset \Omega$ with $\ime(B_0) >
  0$, a.e.\ $\om \in B_0$ is $B$-transitive. (A proof of this can be
  easily devised, for example, with the aid of Fig.~\ref{rwwe-fig}.)
\end{remark}

\proofof{Corollary \ref{cor-evp}} A bounded measurable \fn\ $\phi:
\Omega^\N \longrightarrow \R$ induces a bounded measurable \fn\ $f:
\ps \longrightarrow \R$ via $f(s,\om) := \phi( (\tau_{X_n(s,\om)}
\om)_{n \in \N} ) = \phi( (\om_n)_{n \in \N} )$. By Theorem
\ref{thm-erg}, the asymptotic Birkhoff average of $f$ is constant
$\im$-almost everywhere, i.e., from Proposition \ref{prop-iso}, for
$P_\om$-a.e.\ realization $(x_{n,\om})$ of the RW, in $\ime$-a.e.\
\en\ $\om$. By the definition of the EVP process (compare
(\ref{evp-d2}) with (\ref{rw-d2})), this means that the asymptotic
Birkhoff average of $\phi$ is constant for $\mathcal{P}_{\ime}$-a.e.\
realization $(\om_n)$ of the process. By density, the result extends
to all $\phi \in L^1(\mathcal{P}_{\ime})$.
\qed

\proofof{Corollary \ref{cor-ball}} Let us first notice that, by
(\ref{l-dr}) and the definition of $D$ from Section \ref{sec-dsy},
\begin{equation}
  \label{erg55}
  \ld (\om) = \sum_{i\in\Z^+} q_i(\om) \, d_i = \int_I ds
  \, D(s,\om).
\end{equation}
We apply Theorem \ref{thm-erg} to the displacement \fn\ $D$;
cf.~(\ref{def-xn}):
\begin{eqnarray}
  \label{erg60}
  \lim_{n \to \infty} \frac{X_n (s,\om)}n &=& \lim_{n \to \infty} \, 
  \frac1n \, \sum_{k=0}^{n-1} D \circ \ma^k(s,\om) \nonumber \\
  &=& \int_\ps \im(ds' d\om') D(s',\om') \nonumber \\
  &=& \int_\Omega \ime(d\om') \, \ld (\om'),
\end{eqnarray}
for $\im$- or $\mu$-a.e.~$(s,\om) \in \ps$, that is, $P_\om$-almost
surely for $\Pi$-a.e.\ $\om$ (Proposition \ref{prop-iso}).
\qed

\section{Doubly stochastic environments}
\label{sec-bis}

In the remainder of this note we apply our results to example of the
\emph{doubly stochastic \en s}. These are defined by the condition
that, $\Pi$-almost surely,
\begin{equation}
  \label{bi-en}
  \forall y \in \Z^d, \qquad \sum_{x \in \Z^d} p_{xy} = 1.
\end{equation}
In this case, (A5) is verified by $\Pi$ itself.  In fact, in the
notation of Section \ref{sec-dsy}, (\ref{bi-en}) implies that, for
$\Pi$-a.e.\ $\om$, 
\begin{equation}
  \label{bis10}
  \sum_{i \in \Z_+} q_i (\tau_{-d_i} \, \om) = \sum_{x \in \Z^d}
  p_{0x}(\tau_{-x} \, \om) = \sum_{x \in \Z^d} p_{-x0}(\om) = 1;
\end{equation}
cf.\ (\ref{qi}) and (\ref{zd-action}). Therefore, using
(\ref{steady20}), the invariance of $\Pi$ w.r.t.\ $(\tau_z)$, and
(\ref{bis10}), we obtain, for a measurable $B \subseteq \Omega$,
\begin{eqnarray}
  \int_\Omega \Pi(d\om') \, \mathcal{K} (\om', B) &=& \sum_{i \in
  \Z_+} \int_{\tau_{-d_i(B)}} \hspace{-4pt} \Pi(d\om') \, q_i(\om')
  \nonumber \\
  &=& \sum_{i \in \Z_+} \int_B \Pi(d\om) \, q_i(\tau_{-d_i} \om) 
  \\ \nonumber
  &=& \Pi(B),
\end{eqnarray}
which means precisely that $\Pi$ is a steady state for the EVP; cf.\
(\ref{steady30}). The other condition in (A5) is trivially verified as
$\mathrm{supp}_\Pi \, \Pi = \Omega$.

Therefore, Proposition \ref{prop-im} and all the results of Section
\ref{sec-erg} hold true if $\im$ and $\ime$ are replaced,
respectively, by $\mu$ and $\Pi$. In particular, Theorem
\ref{thm-main} reads:

\begin{proposition}
  \label{prop-bis}
  For a doubly stochastic RWRE verifying \emph{(A1)-(A4)}, the
  annealed process (with law $\mathbb{P}$ and dynamics $\mathcal{F}$)
  and the EVP (with initial state $\Pi$) are stationary and \erg.
\end{proposition}

Suppose further that the RW has \emph{zero local drift}, i.e., $\ld
\equiv 0$; cf.~(\ref{l-dr}). By the invariance of $\Pi$, this is the
same as:
\begin{equation}
  \label{zero-l-dr}
  \forall x \in \Z^d, \qquad \sum_{y \in \Z^d} p_{xy} \, (y-x) = 0,
\end{equation}
for $\Pi$-a.e.\ $\om$. As it is clear, the above means that $( X_n )$
is a martingale. (Examples of doubly stochastic martingales may be
found, for instance, in Appendix A of \cite{l3}.)

In this case, subject to a natural extra condition on the variance of
the jumps, we can prove the \emph{quenched IP}. We state it in the
form of a theorem as soon as we have established some notation.  Given
$( X_n )$, define the continuous \fn\ $R_n : [0,1] \longrightarrow
\R^d$ via the following: For $k = 0, 1, \ldots, n-1$ and $t \in [k/n,
(k+1)/n]$,
\begin{equation}
  \label{rn}
  R_n(t) := \frac {X_k + (nt-k) (X_{k+1} - X_k)} {\sqrt{n}}.
\end{equation}
(Evidently, the graph of $R_n$ is the polyline joining the points
$(k/n, X_k/\sqrt{n})$, for $k = 0, \ldots, n$.) The above can be
regarded as a stochastic process relative to either $P_\om$ (the
quenched \emph{rescaled \tr y}) or $\mathbb{P}$ (the annealed rescaled
\tr y). Then we have the following:

\begin{theorem}
  \label{thm-qip}
  Assume \emph{(A1)-(A4)}, with $\gamma>2$, \emph{(\ref{bi-en})} and
  \emph{(\ref{zero-l-dr})}. Then, for $\Pi$-a.e.\ $\om \in \Omega$,
  the quenched rescaled \tr y $R_n$, relative to $P_\om$, converges to
  the $d$-dimensional Brownian motion with drift 0 and diffusion
  matrix
  \begin{displaymath}
    C := \int_\Omega \Pi (d\om) \sum_{y \in \Z^d}
    p_{0y}(\om) \, y \otimes y.
  \end{displaymath}
  The convergence is intended in the weak-* sense in $C([0,1])$
  endowed with the sup norm. Furthermore, $C$ is positive definite in
  the directions spanned by $\Gamma$.
\end{theorem}

\begin{corollary}
  \label{cor-aip}
  The annealed rescaled \tr y converges to the same Brownian motion as
  in Theorem \ref{thm-qip}.
\end{corollary}
Using Proposition \ref{prop-bis}, the proof of Theorem \ref{thm-qip}
is a standard verification of the hypotheses of the Lindeberg--Feller
Theorem for martingales. (See \cite[Thm.~7.7.4]{d} for a convenient
one-dimensional version of that theorem; cf.\ also \cite[Chap.~4]{hh}
and \cite[Thm.~2.11]{bi}. The multidimensional version follows via the
Cram\'er--Wold device \cite{d,b1}: see, e.g., \cite[p.~1341]{bp}; cf.\
also \cite[Thm.~3.3.4]{z1} and \cite{la}.)  The final assertion on the
diffusion matrix follows easily from the expression given above: if
$\langle v, Cv \rangle = 0$, then, in a.e.\ \en\ $\om$, $\sum_y p_{0y}
\langle y, v \rangle^2 = 0$, so the \rw\ is confined in the orthogonal
subspace of $v$. By (A4), $v$ is orthogonal to
$\mathrm{span}(\Gamma)$.

\skippar

Theorem \ref{thm-qip} can be compared to the main result of \cite{ko},
which is a CLT for doubly stochastic \en s. The standard hypotheses
there (most notably, uniform ellipticity) are stronger than the
present, which makes Theorem \ref{thm-qip} an improvement in the
martingale case. (However, the CLT of \cite{ko} only requires the
\emph{average} of the local drift to be null: $\Pi(\ld) = 0$, allowing
for a much larger class of \rw s.)

An interesting consequence of Theorem \ref{thm-qip} is the almost sure
recurrence in dimension one and two.

\begin{proposition}
  \label{prop-rec}
  Under the hypotheses of Theorem \ref{thm-qip}, with $d \le 2$, the
  \rw\ is $\mathbb{P}$-almost surely recurrent. Equivalently, for
  $\Pi$-a.e.\ $\om \in \Omega$,
  \begin{displaymath}
    \liminf_{n \to \infty} |X_n| = 0,
  \end{displaymath}
  $P_\om$-almost surely.
\end{proposition}
For this we use a general result of Schmidt \cite{s}, as outlined,
e.g., in \cite[Sect.~IV]{l3}.

\end{document}